\documentstyle[12pt,
amstex,amssymb]{article}

\def\GK{\operatorname{GK}}

\def\ЌЋ„{\operatorname{ЌЋ„}}

\def\PI{\mathop{\hbox{PI}}\nolimits}

\newtheorem{theorem}{Theorem}[section]
\newtheorem{corollary}[theorem]{Corollary}

\newtheorem{proposition}[theorem]{Proposition}
\newtheorem{definitionhead}[theorem]{Definition}
\newenvironment{definition}{\begin{definitionhead}%
\normalshape}{\end{definitionhead}}
\def\mytheoremstyle{\it}

\newcommand\Endproof{\nopagebreak\strut%
\nopagebreak\hfill\nopagebreak$\Box$\medbreak}

\def\demo#1{\trivlist\item[\hskip\labelsep{\bf#1\unskip.}]}
\def\enddemo{\endtrivlist}
\def\Pid{\operatorname{PIdeg}}

\def\GK{\operatorname{GKdim}}

\def\Var{\operatorname{Var}}

\def\PI{\operatorname{\text{\it PI}}}

\def\@biblabel#1{#1.}
\def\l@sect#1#2{\addpenalty{\@secpenalty}% good place for page break
   \addvspace{1.0em plus\p@}%
   \@tempdima 1.5em
   \begingroup
     \parindent \z@ \rightskip \@pnumwidth
     \parfillskip -\@pnumwidth
     \bf
     \leavevmode
      \advance\leftskip\@tempdima
      \hskip -\leftskip
     #1\nobreak\hfil \nobreak\hbox to\@pnumwidth{\hss #2}\par
   \endgroup}

\def\nextsection{\refstepcounter{\c@section}{1}}

\begin{document}
\author{ A.~A.~Chilikov\footnote{MIPT, BIU},~A.~Ya.~Belov\footnote{College of Math. and Stat., Shenzhen University}} 

\title {Normal basises of algebras and Exponential Diophantine equations in rings of positive
characteristic}
\date{}
\maketitle
\medskip

UDC 512.5+511


\begin{abstract}

\medskip

In this paper we discourse basises of representable algebras. This question lead to arithmetic problems.
 We prove  algorithmical solvability of exponential-Diophantine
equations in rings represented by matrices over fields of positive
characteristic. Consider the system of exponential-Diophantine equations
$$ \sum\limits_{i=1}^s P_{ij}(n_1,\dots,n_t)b_{ij0}a_{ij1}^{n_1}b_{ij1}\dots
a_{ijt}^{n_t}b_{ijt}=0 $$
where $b_{ijk},a_{ijk}$ are constants from matrix ring of characteristic
$p$, $n_i$ are indeterminates. For any solution $(n_1,\dots,n_t)$ of the
system we construct a word (over an alphabet containing $p^t$ symbols)
${\overline \alpha_0},\dots,{\overline \alpha_q}$ where ${\overline
\alpha_i}$ is a $t$-tuple $\langle n_1^{(i)},\dots,n_t^{(i)}\rangle$,
$n^{(i)}$ is the $i$-th digit in the $p$-adic representation of $n$. The
main result of this paper is as follows: the set of words corresponding in
this sense to solutions of a system of exponential-Diophantine equations is
a regular language (i.e. recognizable by a finite automaton). There exists
an effective algorithm which calculates this language. This algorithm is
constructed in the paper.
\end{abstract}

The research was supported by Russian Science Foundation, Grant No 17-11-01377.

\section{Introduction}

Systems of exponential Diophantine equations (EDE)
\begin{equation}
  \label{diof_system}
  \sum_{i=1}^{s} P_{ij}(n_1,\ldots,n_t)c_{ij1}^{n_i}\cdots c_{ijt}^{n_t}=0,
\end{equation}
where $P_{ij}$ are some polynomials arise in various areas of moderm
mathematics, and in general case, as J.~Robinson has shown, they are
algorithmically undecidable. Yu.~V.~Matiyasevich has proved algorithmical
undecidability for purely Diophantine equations
$$
  P(n_1,\ldots,n_t)=0.
$$
A number of problems reduces to undecidability of some EDE. However it turns
out that if $c_i$ belong to a field of positive characteristic (and even to
a matrix ring), the problem of finding the set of solutions is
algorithmically decidable. Questions arising in this context occur to be
related with formal languages.

     Investigation of bases of algebras is an inspiration for research of
such equations. Suppose $a_1\prec\ldots\prec a_t$ is an ordered set of
generators for an algebra {\bf A}. The order $\prec$ on this set induces
lexicographical ordering on the set of words in $\{ a_i \}.$ A basis {\bf M}
of {\bf A} as of a vector space is called normal if it is generated by {\it
non-decreasable} (that is, not representable by a linear combination of
lesser words) elements. If {\bf A} is a $PI$-algebra (in particular, if {\bf
A} is representable) then due to Shirshov theorem on height there exist
$h=ht({\bf A})$ and a finite tuple $v_1,\ldots,v_s$ such that {\bf M}
consists of elements of the form
$$
  v_{i_1}^{k_1}\cdots v_{i_t}^{k_t}
$$
where $t\le h.$

     In this connection the question arises on the structure of the set
consisting of degree vectors
$\langle k_1,\ldots,k_t\rangle $, in particular for the representable case.
If the algebra {\bf A} is representable and monomial (that is, defining
relations are of the form $u_j=0$ where $u_j$ are some words) then the
problem of normal basis permits in some sense complete answer. We have the
following

{\bf Theorem (test for representability of a monomial algebra).}

{\it A monomial algebra {\bf A} is representable iff {\bf A} has bounded
height over some finite set of words $v_1,\ldots,v_t,$ the set of defining
relations can be divided into a finite number of series
$  v_1^{k_1}\cdots v_t^{k_t}=0$
where
$$
  \sum\limits_{i} P_{ij}(k_1,\ldots,k_t)c_{ij1}^{k_1}\cdots c_{ijt}^{k_t}=0
$$
and each series corresponds to a specific system of EDE.}

This theorem implies, in particular, existence of representable algebras
whose Hilbert series is transcendental as well as algorithmic
undecidability of isomorphism problem for a pair of subalgebras in a matrix
algebra over a polynomial ring.

Nevertheless for positive characteristic the situation is much simpler.
Although Diophantine problems arise more often in this case, their solution
is simpler. The set of solutions of an EDE admits effective description in
terms of $p$-adic decomposition of indeterminates $n_1,\ldots,n_t.$ Since
values of $P_{ij}(n_1,\ldots,n_t)$ are periodical with period
$p$ in each $n_i,$ it suffices to investigate equations of the form
$$
  \sum_{i=1}^{s} c_1^{n_1}\cdots c_t^{n_t}=0.
$$

Consider some solution of an EDE: $\langle  n_1, \ldots, n_t \rangle .$
To each its component $n_i$ attach its $p$-adic decomposition
$n_i^k\ldots n_i^0.$ Thus to each solution we attach the sequence of tuples
of figures
$\langle n_1^k, \ldots, n_t^k\rangle ,$$\ldots,$ $\langle n_1^0, \ldots,
 n_t^0\rangle .$
Interprete these tuples as letters, and our sequences as words over the
alphabet consisting of tuples. The set of all words corresponding to
solutions of EDE forms a language over a finite alphabet. We are ready to
formulate the main result of this paper.

{\bf Theorem 1.}
{\it The set of words corresponding to a system of EDE is a regular
language. In other words, there exists an oriented graph with arrows marked
by letters corresponding to finite tuples of figures (the number of letters
is $p^t$). Some vertex is declared initial, and some other vertices are
declared final. There exists 1-1 correspondence between solutions of our
system and words which may be read along a path consisting of arrows and
going from the initial vertex to a final one. These paths are allowed to
have arbitrary length, and each vertex (including initial and final ones)
may be passed arbitrarily many times.
}

There exists an effective algorithm for constructing such a graph. Below we
present its description.

\section{Bases of Representable and $PI$-algebras}

The earliest purely combinatorial result of this kind occurred to be
{\bf A.~I.~Shirshov height theorem.} {\it Let $A$ be a finitely generated
$\PI$-algebra. Then there exists a finite set of elements
$Y$ and an integer $H\in {\mathbb N}$ such that $A$ is linearly represented by (that is, is generated by linear combinations of) the set of elements of the form
$$
v_1^{k_1} v_2^{k_2}\ldots v_h^{k_h}\quad \hbox{where $h\le H$, $v_i\in
Y$}.
$$
For $Y$ we may take the set of words of degree $\le m$.
}
Such an $Y$ is called {\it a Shirshov basis} of the algebra $A$.

The above theorem implies positive solution of Kurosh problem and of other
Burnside-type problems for $\PI$-rings. In fact, if $Y$ is a Shirshov basis
consisting of algebraic elements then the algebra $A$ is finite-dimensional.
Thus Shirshov theorem explicitly determines the set of elements whose algebraicity
implies algebraicity of the whole algebra. We also have

\begin{corollary}
 If $A$ is a $\PI$-algebra of degree $m$ and all words in its generators of degree
$\le m$ are algebraic then $A$ is locally finite.

\end{corollary}

Height theorem also implies

\begin{corollary}[Berele]
Let $A$ be a finitely generated $\PI$-algebra. Then $\GK(A)<\infty$.

\end{corollary}

$\GK(A)$ is the {\it Gelfand -- Kirillov dimension of the algebra $A$}, that is,
$$\GK(A)=\lim_{n\to\infty}\ln V_A(n)/\ln(n)$$
where $V_A(n)$ is {\it
the growth function of $A$}, that is, the dimension of the vector space generated by
words of degree $\le n$ in generators of $A$.

To prove the corollary, it suffices to observe that the number of solutions of
inequality $k_1 |v_1|+\cdots+k_h|v_h|\le n$ with $h\le H$ does not exceed
$N^{H}$, and so $\GK(A)\le h(A)$.

Thus we obtain various consequences from Height theorem. A little later we
discuss questions concerning conversion of these implications. To begin with, we
introduce some notions and notation.

The number $m=\deg(A)$ will denote {\it the degree of the algebra}, that is, the
minimal degree of an inequality satisfied by it; $n=\Pid(A)$ is {\it the
complexity} of $A$, that is, the maximal $k$ such that
${\mathbb M}_k$, the algebra of matrices of size $k$, belongs to the variety
$\Var(A)$ generated by $A$.

It is convenient to replace the notion of  {\it height} by a close notion of {\it
essential height}.

\begin{definition}
An algebra $A$ has {\it essential height $h$} over a finite set $Y$ which is called
an {\it $s$-basis} if there exists a finite set $D\subset A$ such that $A$ is
linearly representable by elements of the form $t_1\cdot\ldots\cdot t_l$ where $l\le
2h+1$ and $\forall i (t_i\!\in\! D \vee t_i=y_i^{k_i};y_i\in Y)$, and the set of
$i$ having $t_i\not\in D$ contains $\le h$ elements.
\end{definition}

Loosely speaking, each long word is a product of periodical parts and of ``layers''
having bounded length. Essential height is the number of these periodical pieces, and ordinary height depends also upon ``layers''.

Height theorem gives rise to following questions:

\begin{enumerate}

\item To what classes of rings Height theorem may be extended?

\item For which $Y$ the algebra $A$ has bounded height?

Since now, we consider the associative case.

\item How to evaluate height?

\item What does the degree vector $(k_1,\ldots,k_h)$ look like?
First of all, which sets of its components are essential, that is, which sets of
$k_i$ can be simultaneously unbounded? What is the value of essential height?

\item A question regarding finer structure of the set of degree vectors: does it
have any regularity properties?

At last, the range of questions forming the subject of this paper.

\item Which sets of words can be chosen for $\{v_i\}$?

\end{enumerate}

Now we proceed to discuss the above questions.
%\medskip

{\bf Non-associative generalizations.} Height theorem has been extended to
certain classes of rings close to associative rings. S.~V.~Pchelintsev
has proved it for alternative and $(-1,1)$ cases,
S.~P.~Mishchenko  has obtained an analogue of Height theorem
for Lie algebras with a sparse identity. The author  has proved Height
theorem for a certain class of rings asymptotically close to associative rings
and in particular including alternative and Jordan $\PI$-algebras.

{\bf Shirshov bases.} Let $A$ be a $\PI$-algebra, and suppose a subset
$M\subseteq A$ is its $s$-basis. Then if all elements of
$M$ are algebraic over $K$ then $A$ is finite-dimensional
(Kurosh problem). Boundedness of essential height over $Y$
implies ``positive solution of Kurosh problem over $Y$''. The converse is much less
trivial.

\begin{theorem}  [A.~Ya.~Belov]         \label{ThKurHmg}
Suppose $A$ is a graded $\PI$-algebra, $Y$ is a finite set of homogeneous elements.
Then if the algebra $A/Y^{(n)}$ is nilpotent for each $n$ then
$Y$ is an $s$-basis for $A$. If in this situation $Y$ generates
$A$ as an algebra then $Y$ is a Shirshov basis for $A$.

\end{theorem}

($Y^{(n)}$ denotes the ideal generated by $n$th powers of elements from $Y$.)

The following example demonstrates that the straightforward converse of Kurosh
problem for non-graded case does not have positive solution. Suppose
$A={\mathbb Q}[x,1/x]$. Each projection $\pi$ such that $\pi(x)$ is algebraic has
finite-dimensional image. Nevertheless the set $\{x\}$ is not an $s$-basis for the
algebra ${\mathbb Q}[x,1/x]$.

Thus the definition of Kurosh set is chosen as follows:

\begin{definition}
A set $M\subset A$ is called {\it a Kurosh set} if each projection
$\pi\colon A\otimes K[X]\to A'$ having image $\pi(M)$ integral
over $\pi(K[X])$ is finite-dimensional over $\pi(K[X])$.
\end{definition}

We proceed to formulate a generalization of this theorem for non-homogeneous case.

\begin{theorem}[A.~Ya.~Belov]           \label{ThKurGen}
Let $A$ be a $\PI$-algebra, $M\subseteq A$ a Kurosh subset in $A$. Then $M$ is an
$s$-basis for $A$.
\end{theorem}

The following proposition shows that Theorem \ref{ThKurGen} is a generalization of
Theorem \ref{ThKurHmg}:

\begin{proposition}
Let $A$ be a graded algebra, $Y$ a set of homogeneous elements. Then if the algebra
$A/Y^{(n)}$ is locally nilpotent for all $n$ then $Y$ is a Kurosh set.

\end{proposition}

Thus boundedness of essential height is a non-commutative generalization of
{\it integrity}.

\demo{Remarks}
a) Note that in the case of Lie $\PI$-algebras, Kurosh problem has positive solution
but Height theorem fails.

b) The theorem extends to some class of rings asymptotically close to associative
rings (with bounded $l$-length, finitely generated algebra of left multiplications,
and associative powers).
\enddemo

{\bf Estimates of height.} The original A.~I.~Shirshov's proof was purely
combinatorial (based on elimination technique developed by him for Lie algebras, in
particular in the proof of Freedom theorem), however it did not provide any
reasonable estimates for height. Later A.~T.~Kolotov
obtained an estimate for $ht(A)\le s^{s^m}$\
($m=\deg(A)$,\, $s$ is the number of generators). Subsequently, E.~I.~Zel'manov
\cite{Dnestrovsk} raised the question on existing of an exponential estimate which
was obtained later on by the Belov.

\demo{Shirshov height theorem}                 %\label{Th1b28}%
{\mytheoremstyle
Suppose $A$ is an $l$-generated $\PI$-algebra of degree $m$. Then the height of $A$
over the set of words having degree $\le m$ is bounded by a function $H(m,l)$
where $H(m,l) < 2ml^{m+1}$.
}
\enddemo

{\bf Essential height.} Clearly, essential heght is an estimate for Gelfand --
Kirillov dimension and an $s$-basis is a Shirshov basis iff it generates  $A$ as an
algebra.

In the representable case the converse is true.

\begin{theorem}[A.~Ya.~Belov \cite{BBL}]
Suppose $A$ is a finitely generated representable algebra and
$H_{Ess}{}_Y(A)<\infty$. Then $H_{Ess}{}_Y(A)=\GK(A)$.

\end{theorem}

\begin{corollary}[V.~T.~Markov]
The Gelfand -- Kirillov dimension of a finitely generated representable algebra is an integer.

\end{corollary}

\begin{corollary}
If $H_{Ess}{}_Y(A)<\infty$ amd an algebra $A$ is representable then
$H_{Ess}{}_Y(A)$ is independent of the $s$-basis $Y$.

\end{corollary}

Due to local representability of relatively free algebras, the Gelfand -- Kirillov
dimension in this case also equals the essential height.

{\bf Structure of degree vectors.} Thus in the representable case both Gelfand --
Kirillov dimension and essential height behave well. Nevertheless even in this case
the set of degree vectors can have bad structure, namely, it can be the complement
for the set of solutions for some system of exponential-polynomial Diophantine
equations. Consequently, there exists an example of a representable algebra having
transcendent Hilbert series. However in the case of relatively free algebra the
Hilbert series is rational.

{\bf Shirshov bases consisting of  words.} Their description is given by the
following theorem:

\begin{theorem}[A.~Ya.~Belov]            \label{ThBelheight}
A set of words $Y$ is a Shirshov basis of an algebra $A$ iff for each word
$u$ having length $\le m = \Pid(A)$, the complexity of $A$, the set $Y$ contains some word which is cyclically conjugate to some degree of $u$.

\end{theorem}

A.~I.~Shirshov himself has shown that for a Shirshov basis we may take the set of
words having
degree at most $\deg(A)$. I.~V.~Lvov has proved boundedness of height over the
set of words having length at most $\deg(A)-1$.
S.~Amitsur and I.~P.~Shestakov conjectured that if all words having length not
exceeding the complexity $\Pid(A)$ are algebraic then the algebra is
finite-dimensional.
I.~V.~Lvov reduced this statement to the following:

\begin{theorem}
Let $A$ be a finite-dimensional subalgebra in the matrix algebra of order
$n$, and let $a_1,\dots,a_s$ be its generators. Then if all words in $a_1,\dots,a_s$
having degree $\le n$ are nilpotent then $A$ itself is nilpotent.

\end{theorem}

Note that $n$ is the precise estimate.

Shestakov's conjecture was proved by V.~A.~Ufnarovsky  and
by G.~P.~Chekanu.
\footnote{From a private letter by the latter: ``we worked in the same area $\dots$
We both have stood this friendly and creative concurrence (we did begin this
deliberately, with agreement that we work in different languages)''. The proofs were
based on ``the spirit of independence''. Subsequent papers of these authors contained various specifications and generalizations of these theorems \cite{Chekanu3}.}.
Later the author  \cite{Belov1} showed that for $\{v_i\}$, we may take the set of
words from Shestakov's conjecture. This result also was announced by G.~P.~Chekanu. Later on, another proof of this fact was obtained by
V.~Drensky.

In the sequel, we focus on the range of problems concerned to relations between
Height theorem and Independence theorem.

Independence theorem may be formulated, in particular, as follows
 %\cite{Chekanu}, \cite{Uf7}:

\begin{theorem}[Independence theorem]
Suppose the following is true:
\begin{enumerate}
\item a word $W = a_{i_{1}}\ldots a_{i_{n}}$
is the minimal word in the left lexicographical ordering on the set of all nonzero
products
having length $\le n$;
\item the extreme parts of $W$ are nilpotent.
\end{enumerate}
Then initial subwords of $W$ are linearly independent.

\end{theorem}

To deduce I.~P.~Shestakov's conjecture (or, equivalently, I.~V.~L'vov's statement)
from this theorem, it suffices to consider a faithful representation of $A$ by
operators on $n$-dimensional space $V$. Let $v_{1},\ldots,v_{n}$ be a basis of this
space, then for some $v_i$ we have $m_iW\ne 0$. Consider the auxiliary algebra
generated by
$V$ and $A$. Suppose $V\cdot V=A\cdot V=0$ and the action of $VA$ coincides with
module
multiplication. Reorder the generators as follows: $v_1\succ\dots\succ v_n\succ
a_1\succ\dots\succ a_s$ and apply Independence theorem.
\Endproof

Original proofs of Independence theorem were rather complicated. Application of
symbolic dynamics technique involving infinite words or {\it superwords} allowed to
clarify them. Technique of superwords occurred to be rather close to the lines of
structure theory. Its role does not reduce to proving statements like
Independence
theorem. Using superwords allows to prove Height theorem, nilpotence of the Lie
algebra generated by sandwiches \cite{Uf7}, coincidence of nilradical and Jacobson
radical in monomial algebras, to describe bases of algebras with extremal growth
function $V_A(n)=\frac{n(n+3)}{2}$, and also to describe weakly Noetherian,
semisimple and semiprimary monomial algebras \cite{BBL} and to obtain some other
combinatorial results in the theories of semigroups and rings.

Many properties of algebras are defined by monomial relations. For example, such are
the conditions of Shestakov's conjecture, namely, nilpotence of  words whose degree
does not exceed complexity.

This conjecture is related to the structure of the matrix algebra. Multiplication of
matrix units $E_{ij}$ is almost monomial, and the language of representations of
matrix algebras clarifies ``matrix'' properties of semisimple components. It is no
coincidence that many authors dealing with independence actively used a similar
technique of matrix constructions for other problems concerning local finiteness
\cite{Chekanu4}.

\section{Preliminries}

Recall now some facts from the theory of formal languages.

{\bf Definitions.}
{\it A finite automaton} is a finite oriented graph some vertex of which is
declared initial, some vertices are declared final, and each edge is marked
by a symbol of some finite alphabet.

{\it A regular language} is a set of words which may be read at edges of
some finite automaton along a path from the initial vertex to a final one.
We say that this automaton {\it represents} the given language.

{\it A concatenation} $vu$ of two words $u$ and $v$ is obtained by adding
$v$ to $u$ (in our case from the left).
{\it A concatenation} of two languages $L_1$ and $L_2$ is the language $L=$
\{$uv \mid  u\in L_1, v\in L_2$\}.

{\it The closure} $L^*$ of a language $L$ is the set of all powers of words
from $L.$

{\it An atomary language} consists of a single word consisting of a single
letter.

One of the simplest instances of regular languages is the set of all words
not including subwords from a fixed finite list. A description of regular
languages in terms of operations over languages is provided by the following

{\bf Theorem (Cleenee).}
{\it A language is regular iff it can be obtained from atomary languages by
finite number of operations of joint, meet, comcatenation and closure.}

For more detail and for the proof of Cleenee's theorem see Salomaa
[1, p.~24-37].
\section{Basic notation and constructions}

In the sequel, we use following notation.\\
   $\overline\vartheta=(\vartheta_1,\ldots,\vartheta_r)$ is a tuple of
variables.\\
   $\overline n-(n_1,\ldots,n_t)$ is a tuple of indeterminates.\\
   {\bf N} is the set of natural numbers.\\
   $\Sigma_1$ is the alphabet consisting of tuples of figures
$0,1,\dots,p-1$ having length $t$.\\
   $\Sigma_2$ is the alphabet consisting of tuples of figures
$0,1,\dots,p-1$ having length $r$.\\
   $\Sigma_1^*$ is the set of all finite words over $\Sigma_1.$\\
   $\Sigma_2^*$ is the set of all finite words over $\Sigma_2.$\\
   Words from $\Sigma_1^*,\Sigma_2^*$ are written from the right to the
left.\\
   $\lambda$ is the empty word.\\
   $l(u)$ is the length of the word $u.$\\
   ${\bf R}={\bf Z_p}[\vartheta_1,\ldots,\vartheta_r]$ is the ring of
polynomials over ${\bf Z_p}.$\\
   {\bf F} is the quotient field for {\bf R}.\\
   {\bf A} is the algebraic closure for {\bf F}.\\
   Since to each sequence of figures with radix $p$ there corresponds a
number from {\bf N}, the following maps are well-defined:\\
   $\phi:\Sigma_1^*\rightarrow N^t$ which maps any word from
$\Sigma_1^*$ to a tuple consisting of $t$ numbers written with radix $p$,\\
   $\psi:\Sigma_2^*\rightarrow N^r$ which maps any word from
$\Sigma_2^*$ to a tuple consisting of $r$ numbers written with radix $p$.\\
   $\phi^{(i)}:\Sigma_1^*\rightarrow {\bf N}$ is $i$th component of
$\phi.$\\
   $\psi^{(i)}:\Sigma_2^*\rightarrow {\bf N}$ is $i$th component of
$\psi.$\\
   $\overline f=(f_1,\ldots,f_t)$ is a tuple of polynomials.\\
   In the sequel, we also denote by ${\overline f}^p$ the tuple consisting
of $p$th powers of $f_i$: $(f_1^p,\ldots,f_t^p),$ and we denote by
${\overline\vartheta}^p$ the tuple consisting of $p$th powers of
$\vartheta_i$:
$(\vartheta_1^p,\ldots,\vartheta_r^p).$\\
   Products of the form $f_1^{\phi^{(1)}(u)}\cdots f_t^{\phi^{(t)}(u)}$ and
$\vartheta_1^{\psi^{(1)}(v)}\cdots \vartheta_r^{\psi^{(r)}(v)}$ will be
denoted ${\overline f}^{\phi(u)}$ and ${\overline\vartheta}^{\psi(v)}$
accordingly.\\
\section{ Equations over a ring of polynomials}

Let {\bf K} be a matrix ring having positive characteristic $p$.
We prove now an auxiliary statement which allows to reduce the class of
considered equations.

{\bf Proposition.}
{\it If for each equation of the form
\begin{equation}
  \label{ur_wida}
  \sum_{i=1}^{s} b_{i0} a_{i1}^{n_1} b_{i1} \cdots a_{it}^{n_t} b_{it} =0
\end{equation}
having coefficients from {\bf K} the set of words corresponding to its
solutions (as it was defined in Introduction) is a regular language then the set of words corresponding to solutions of any system of equations of the
form
\begin{equation}
  \label{ur_wida2}
   \sum_{i=1}^{s} P_{ij}(n_1,\ldots,n_t) b_{ij0} a_{ij1}^{n_1} b_{ij1}
\cdots a_{ijt}^{n_t} b_{ijt} =0
\end{equation}
over {\bf K} is a regular language.}

{\bf Proof.}
   First note that the set of solutions of a system is the meet of sets of
solutions for equations of the system. So by virtue of Cleenee's theorem,
regularity of languages corresponding to single equations of the form
(\ref{ur_wida2}) implies regularity of languages corresponding to systems of
such equations.

     Consider now an equation of the form (\ref{ur_wida2}). Let
$\langle n_1,\ldots,n_t\rangle $ be a tuple of numbers
$n_i=n_i^0+p n_i^{'}$ where $n_i^0$ is the last digit with radix $p$ in the
number $n_i.$ For fixed $\langle n_1^{'},\ldots,n_t^{'}\rangle $ we have
$$
P(n_1^{'},\ldots,n_t^{'})=P(n_1,\ldots,n_t),
$$
so
$\langle n_1,\ldots,n_t\rangle $ is a solution of an equation of the form
(\ref{ur_wida2}) iff $\langle n_1^{'},\ldots,n_t^{'}\rangle $ is a solution
of an equation of the form (\ref{ur_wida}). Regularity of the set of all
$\langle n_1^{'},\ldots,n_t^{'}\rangle $ obviously implies regularity of the
set of all $\langle n_1^{'},\ldots,n_t^{'}\rangle $ since the corresponding
words are obtained by adding the tuple
$\langle n_1^0,\ldots,n_t^0\rangle .$
Finally, the complete set of solutions of the original equation is the joint
of sets of solutions corresponding to distinct tuples
$\langle n_1^0,\ldots,n_t^0\rangle .$
Again by Cleene's theorem we have regularity of languages corresponding to
any equations of the form (\ref{ur_wida2}) over {\bf K}. The proof is
complete.

Thus we have reduced investigation of solutions for some system of EDE to
the case of a single equation which furthermore has no polynomial (in $n$)
parts.

Consider an EDE over {\bf R}:
\begin{equation}
  \label{edu}
  \sum_{i=1}^{s} Q_i(\overline\vartheta)[P_{i1}^{n_1}](\overline\vartheta)
\cdots[P_{it}^{n_t}](\overline\vartheta)=0.
\end{equation}
Its solution is a tuple of numbers  $\overline n=(n_1,\ldots,n_t),
n_i\in{\bf N}.$

{\bf Definition.}
{\it A word-solution} of the EDE (\ref{edu}) is $u\in\Sigma_1^*$ such that
$\phi(u)=\overline{n}$ where $\overline{n}$ is a solution for (\ref{edu}).

Now we may write the equation in $u$
\begin{equation}
  \label{ur_otn_u}
  \sum_{i=1}^{s} Q_i(\overline\vartheta)
[{\overline{P_i}}^{\phi(u)}](\overline\vartheta)=0.
\end{equation}
     In the sequel, word-solutions will also be called solutions simply.
The main result of this Section may be stated in the form of the following

{\bf Theorem 1.}
{\it
 $L\subset\Sigma_1^*$, the set of word-solutions for the equation
(\ref{ur_otn_u}), is a regular language.}

      Before presenting the proof, we describe its basic idea. Let
$Q(x)$ be a polynomial having coefficients from
${\bf Z_p}.$ Let us investigate the result of removing brackets in
$Q^n(x).$ Write $n$ with radix $p$:
$$
 n=n_{0}+n_{1}p+\ldots+n_{k}p^{k}+n_{k+1}p^{k+1}+\ldots+n_{s}p^{s}.
$$
 Put $Q_{k}(x)=Q^{k}$ where $k=0,1,\ldots,p-1.$ Clearly $Q_{0}=1,$
$Q_{1}=Q.$ Then
$$
 Q(x)^n=Q_{n_{0}}(x)Q_{n_{1}}(x^p)\cdots Q_{n_{k}}(x^{p^{k}})
 Q_{n_{k+1}}(x^{p^{k+1}})\cdots Q_{n_{s}}(x^{p^{s}})~~~(*)
$$
 since ${Q(x)}^{p^k}=Q(x^{p^k}).$
 Consider a section of the product
$$
 R_k=Q_{n_0}(x)\cdots Q_{n_k}(x^{p^k}).~~~(**)
$$
Collect terms with the same remainder $\alpha$ of the degree of
$x$ modulo $p^{k+1}$. In other words, represent the section as sum
$$
 \sum_{\alpha=0}^{p^{k+1}-1} x^{\alpha}R_{\alpha}(x^{p^{k+1}}).~~~(***)
$$
Now note the following.
\\
1. Multiplication by the rest of the product $(*)$ (that is, by
$Q_{n_{k+1}}(x^{p^{k+1}})\cdots Q_{n_{s}}(x^{p^{s}})$) does not lead to
cancellation of terms in $(***)$ having distinct $\alpha$.\\
2. The degree of $R_k$ does not exceed $ (deg \mbox{ } Q)(p-1)p^k.$
Hence the degree of $R_{\alpha}$ does not exceed $ (deg \mbox{ } Q)(p-1).$
Furthermore, since we work over a finite field, the number of distinct types
of $R_{\alpha}$ ({\it small types}) is bounded (and does not exceed
$p^{(deg \mbox{ } Q)(p-1)}$).\\
3. Due to observation 1, we need not all the information on the sum $(***)$
but only the following: which polynomials $R_{\alpha}$ do exist for given
$k.$
\\
The set of all existing polynomials will be called
{\it a large type.} Clearly the number of large types is finite
(and does not exceed $2^{p^{ (deg \mbox{ } Q)(p-1)}}$).
It is also clear that the large type for a given $k$ uniquely determines the large type for $k+1.$ This implies finiteness of the space of large types for products $(**).$ If we use several monomials then we have to define
{\it the small type} of the sum
 $$
  S_1 Q_1^{n_1}+\ldots+ S_l Q_l^{n_l}
 $$
as the tuple consisting of small types of summands, and
{\it the large type} as the tuple consisting of involved small types.

     Now if we consider polynomials in several variables and products of the
form
 $$
  S_{i} Q_{i1}^{n_1}\cdots Q_{it}^{n_t}
 $$
then we have to take tuples of remainders modulo $p^{k+1}$ and to collect
variables having corresponding degrees. In this case we have:
 $$
  \sum_{\alpha_1,\ldots,\alpha_r} x_1^{\alpha_1}\cdots x_r^{\alpha_r}
  R_{\overline\alpha} (x_1^{p^{k+1}},\ldots,x_r^{p^{k+1}})
 $$
for each monomial. Then the  small type of the monomial is
$R_{\overline\alpha}$, the small type of the system is the tuple consisting
of small types of monomials with given $\alpha,$ and finally the large type
of the system is the tuple consisting of involved small types.
It is easily seen that writing a new figure to the left from variables
$n_i$ corresponds to change of large types (depending of the written
figure), and vanishing of the expression
 $$
  \sum R_i Q_{i1}^{n_1}\cdots Q_{it}^{n_t}
 $$
depends only on its large types (since vanishing of its components
obtained by grouping terms described above depends only on its small
types). Thus we obtain a finite graph of states. Its vertices are large
types, and arrows marked by tuples of figures $0,\dots,p-1$ are
transformations of large types. The initial vertex corresponds to the large
type of zero, and final vertices correspond to those large types which
provide cancellation of all summands.

     Clearly there is a correspondence between words which may be read at
arrows of the graph along a path from the initial vertex to a final, and
solutions of the EDE
$$
 \sum R_i Q_{i1}^{n_1}\cdots Q_{it}^{n_t}=0.
$$

     Now we turn to formal details. First we introduce some important
constructions. Let $f$ be a polynomial from {\bf R}. Then
$$
f(\overline\vartheta)=\sum\limits_{y\in\Sigma_2}
f_y({\overline\vartheta}^p){\overline\vartheta}^{\psi(y)},
$$
and this decomposition is unique.

{\bf Definition.}
{\it The weeding} of a polynomial $f$ by a symbol $y$ is the polynomial
$\varepsilon _y(f)=f_y(\overline\vartheta).$
{\it The weeding} of a polynomial $f$ by a word $v=y_k\cdots y_0$ is the
polynomial
$\varepsilon_v(f)=\varepsilon_{y_k}(\cdots(\varepsilon_{y_0}(f))\cdots).$

{\bf Remark.}
Weeding is a way to collect, as it was mentioned above, polynomials in which
degrees of variables coincide modulo $p^k$.

It is easy to see that
 $ \varepsilon_y(f+g)=\varepsilon_y(f)+\varepsilon_y(g)$
 and $deg \mbox{ }\varepsilon_y(f)\le\frac{1}{p}  deg \mbox{ } f.$

{\bf Lemma 1.}
{\it
 $\varepsilon_y(f(\overline\vartheta)g({\overline\vartheta}^p))=
 \varepsilon_y(f(\overline\vartheta))g(\overline\vartheta).$\\
In other words,in weeding polynomials in ${\overline\vartheta}^p$  are
factored out loosing degree $p.$}

{\bf Proof.}
Represent $f$ in the form
$\sum\limits_{y\in\Sigma_2}{\overline\vartheta}^{\psi(y)}
f_y({\overline\vartheta}^p).$ Then
$$
 f(\overline\vartheta)g({\overline\vartheta}^p)=
 \sum\limits_{y\in\Sigma_2}{\overline\vartheta}^{\psi(y)}f_y
({\overline\vartheta}^p)g({\overline\vartheta}^p)=
 \sum\limits_{y\in\Sigma_2}{\overline\vartheta}^{\psi(y)}(f_y({\overline
\vartheta}^p)g({\overline\vartheta}^p)).
$$
By definition of weeding we immediately obtain
$$
 \varepsilon_y(f(\overline\vartheta)g({\overline\vartheta}^p))
 =\varepsilon_y(f(\overline\vartheta))g(\overline\vartheta).
$$

{\bf Lemma 2.}
{\it
 Let $c$ be a positive integer. Then a polynomial $f(\overline\vartheta)$
vanishes iff every its weeding by a word of length $c$ vanishes.}

{\bf Proof.}
Use induction in $c.$

 a) Base of induction. Suppose $c=1.$ Then $f=0$ obviously implies
$\varepsilon_y(f)=0.$ \\
    Conversely, suppose $\varepsilon_y(f)=0$ for any $ y\in\Sigma_2$.
     Then $ f=\sum\limits_{y\in\Sigma_2}
     \varepsilon_y(f){\overline\vartheta}^{\psi(y)}=0.$
     The base of induction is proved.

 b) The inductive step. Suppose the statement is valid for $c=k$. Then $f=0$
means that
 $\varepsilon_y(f)=0$ for any $y\in\Sigma_2.$ This in turn is equivalent to
$\varepsilon_y(\varepsilon_{v'}(f))=0$ for any
 $ y\in\Sigma_2$ and any $ v'\in\Sigma_2^*$ of length $k.$ Hence
 $\varepsilon_v(f)=0$ for any $ v\in\Sigma_2^*$ of length $k+1.$ The
inductive step is proved.
\\

{\bf Special operators of an equation.}
\\

\nobreak
 Return to the original equation (\ref{ur_otn_u})
$$
 \sum_{i=1}^{s} Q_i(\overline\vartheta)
[{\overline{P_i}}^{\phi(u)}](\overline\vartheta)=0.
$$

{\bf Definition.}
 Suppose $i$ is an integer in the interval from $1$ to $s,$ $u$ and $v$ are
two words of equal length over alphabets $\Sigma_1$ and $\Sigma_2$
accordingly. Then define the {\it special operator} of the equation
(\ref{ur_otn_u}) $S_{u,v}^{(i)}(f)$
 as $\varepsilon_v(f[{\overline{P_i}}^{\phi(u)}]).$

{\bf Lemma 3.}
{\it
 If length of words $u_1,$ $v_1$ is equal, and similarly for $u_2,$ $v_2,$
then
 $$
   S_{u_1 u_2,v_1 v_2}^{(i)}(f)=S_{u_1,v_1}^{(i)}S_{u_2,v_2}^{(i)}(f),
 $$
that is, the special operator corresponding to concatenation is the
composition of special operators corresponding to its factors.}

{\bf Proof.}
Denote by $k$ the length of $u_1$ (equal to the  length of $v_1$). Then the
composition
$S_{u_1,v_1}^{(i)}S_{u_2,v_2}^{(i)}(f)$ equals
  $\varepsilon_{v_1}(\varepsilon_{v_2}( f[{\overline{P_i}}^{\phi(u_2)}] )
  [{\overline{P_i}}^{\phi(u_1)}] ).$
Include $[{\overline{P_i}}^{\phi(u_1)}]$ in the weeding. We have

  $$
  \varepsilon_{v_1}(\varepsilon_{v_2}(
f[{\overline{P_i}}^{\phi(u_2)+\phi(u_1)p^k}]))
  =\varepsilon_{v_1 v_2}( f[{\overline{P_i}}^{\phi(u_1 u_2)}] )=S_{u_1
u_2,v_1 v_2}^{(i)}(f).
  $$
Lemma is proved.

{\bf Lemma 4 (on decreasing the degree).}
{\it
There exists $N_0$ such that for any $N'\ge N_0, 1\le i\le s,$ and for any
 $u\in\Sigma_1^*,v\in\Sigma_2^*$ of equal length $ deg \mbox{ } f\le N'$
 implies $ deg \mbox{ } S_{u,v}^{(i)}(f)\le N'.$\\
In other words, rather large degrees can only decrease under the action of
special operators.}

{\bf Proof.}
The idea of proof is as follows: in a special operator, multiplying by fixed
polynomials increases the degree of the original polynomial not more by a
constant, and after that weeding decreases its degree not less than $p$
times. We proceed to formalize this argument.

     Denote $max \mbox{ } deg \mbox{ } P_{ik}$ by $M.$ Then the required
$N_0$ equals $\frac{prM}{p-1}.$ Indeed, $N'=N_0+K.$ Then if
$ deg \mbox{ } f\le N',$ then for all $x\in\Sigma_1$ we have
       $$
         deg \mbox{ } {\overline{P_i}}^{\phi(x)}=
         \sum\limits_{k=1}^r ( deg \mbox{ } P_{ik})\phi^{(k)}(x)
         \le\sum\limits_{k=1}^{r}Mp=Mpr.
       $$
       Then \\
       $ deg \mbox{ } S_{x,y}^{(i)}(f)= deg \mbox{ }
\varepsilon_y(f{\overline{P_i}}^{\phi(x)})
       \le\left(\frac{deg \mbox{ } f+ prM}{p}\right)$
       $\le rM+\frac{N_0+K}{p}=\frac{prM}{p-1}+\frac{K}{p}\le N'.$\\
The assertion of Lemma now follows.
\\

{\bf Types and their extensions.}
\\

\nobreak
{\bf Definitions.}
 {\it A small type} $T=(f_1,\ldots,f_s)$ is a string of polynomials from
{\bf R} having degree not exceeding
 $N_1= max\{ max$ $deg \mbox{ } (Q_i),N_0\}.$\\
 {\it A large type} T is an arbitrary set of small types.\\
 {\it The extension}  $\pi(u,v)\tau$ of a small type
$\tau=(f_1,\ldots,f_s)$
 by a pair of words $u\in\Sigma_1^*,\mbox{ } v\in\Sigma_2^*$ having equal
length is the small type $\tau'=(f'_1,\ldots,f'_s)$ where
$f'_i=S_{u,v}^{(i)}(f_i).$\\
 {\it The extension} $\Pi (u) T$ of a large type T by a word
$u\in\Sigma_1^*$
is the large type $T'=\{\pi(u,w)\tau \mid
 \tau\in T,w\in\Sigma_2^*,$ $l(w)=l(u)\}.$\\

{\bf Remark.}
It is easy to observe that the operation of extension is defined for all
small types. Indeed, if $ deg \mbox{ } f_1\le N_1$ where $N_1\ge N_0$ then
$deg \mbox{ } f_i'\le N_1,$ so $\tau'$ is also a small type.

 Moreover small types are strings of polynomials of bounded degree in
$r$ variables over a finite field, so their number is finite. The number of
large types is finite as well since they are subsets of a finite set. \\

{\bf Definitions.}
 {\it The small type of a pair of words} $u\in\Sigma_1^*, \mbox{ }
v\in\Sigma_2^*$ of equal length
 is $\tau(u,v)=\pi(u,v)\tau(\lambda,\lambda)$
 where $\tau(\lambda,\lambda)=(Q_1,\ldots,Q_s).$\\
 Define {\it the large type of a word} $u\in\Sigma_1^*$ as $T(u)=
 \{\tau(u,w) ; l(w)=l(u)\}.$

{\bf Lemma 5.}
{\it
$T(u)=\Pi(u)T(\lambda),$
that is, the large type of a word $u$ may be obtained as an extension by $u$
of the large type of the empty word.}

{\bf Proof.}
By definition $T(\lambda)=\{\tau(\lambda,\lambda)\}.$
 Denote by $T'$ the extension of the type $T(\lambda)$ by $u.$
 Then $T'$ is the set of various extensions
 $\pi(u,v)\tau(\lambda,\lambda)$ of the type $\tau(\lambda,\lambda)$
 by pairs of words $u,v$ where $v$ is an arbitrary word of the same length
as $u.$
Since the small type $\tau(u,v)$ is by definition
 $\pi(u,v)\tau(\lambda,\lambda)$ then $T'$ is the set of all $\tau(u,v)$
where $v$ is an arbitrary word of the same length as $u,$ and so it
coincides with $T(u).$ Lemma is proved.

{\bf Definitions.}
 A small type $\tau=(f_1,\ldots,f_s)$ is {\it good}
 if $\sum\limits_{i=1}^s f_i =0.$\\
 A large type T is {\it good} if all $\tau\in T$ are good.
\\

We proceed to prove the following

{\bf Theorem 2.}
{\it
 A large type $T(u)$ is good iff $u$ is a solution of the equation
 (\ref{ur_otn_u}).}

{\bf Proof.}
Denote the length of $u$ by $c.$ A large type $T(u)$ is good iff
for all $v\in\Sigma_2^*$ having length $c$ the small type
 $\pi(u,v)\tau(\lambda,\lambda)$ is good. This in turn means that for all
such $v$ we have
 $$
 \sum\limits_{i=1}^{s} S_{u,v}^i( Q_i(\overline\vartheta) )=0
 $$
or equivalently
 $$
 \sum\limits_{i=1}^{s} \varepsilon_v( Q_i[{\overline{P_i}}^{\phi(u)}] )=0.
 $$
Furthermore due to linearity of weeding we have
 $$
 \varepsilon_v( \sum\limits_{i=1}^{s} Q_i[{\overline{P_i}}^{\phi(u)}] )=0
 $$
and by lemma 2
 $$
  \sum\limits_{i=1}^{s} Q_i(\overline\vartheta)
  [{\overline{P_i}}^{\phi(u)}](\overline\vartheta)=0.
 $$
But this means that $u$ is a solution of (\ref{ur_otn_u}).

{\bf Lemma 6.}
{\it
  a) $\pi(u_1 u_2,v_1 v_2)\tau=\pi(u_1,v_1)\pi(u_2,v_2)\tau.$\\
 b) $\Pi(u_1 u_2)T=\Pi(u_1)\Pi(u_2)T.$

In other words, an extension of a type by a concatenation is a composition
of extensions by factors.}

{\bf Proof.}
 a) Suppose $\tau=(f_1,\ldots,f_s).$
Then let $\tau'=(f'_1,\ldots,f'_s)$ denote the extension of
     $\tau$ by a pair of words $u_2,v_2$, and let
$\tau''=(f''_1,\ldots,f''_s)$ denote the extension of
     $\tau'$ by a pair of words $u_1,v_1.$ We proceed to prove that
     $\pi(u_1 u_2,v_1 v_2)\tau=\tau''.$ Note that
     $f'_i=S_{u_2,v_2}^{(i)}(f_i)$ and $f''_i=S_{u_1,v_1}{(i)}(f'_i).$
     Hence $f''_i=S_{u_1,v_1}^{(i)}S_{u_2,v_2}^{(i)}(f_i)=
     S_{u_1 u_2,v_1 v_2}^{(i)}(f_i).$\\
     Thus $\pi(u_1 u_2,v_1 v_2)\tau=\tau''=\pi(u_1,v_1)\pi(u_2,v_2)\tau.$
     First assertion of Lemma is proved.\\
 b) Denote $\Pi(u_2)T$ by $T',$ and $\Pi(u_1)T'$ by $T''.$
     We shall prove that $\Pi(u_1)\Pi(u_2)T=T''=\Pi(u_1 u_2)T.$
By definition, $T'$ consists of various extensions of types from
     T by pairs of words $u_2,v_2$ having equal length. Also by definition,
     $T''$ is the set of extensions of types from $T'$ by pairs of words
$u_1,v_1$ having equal length. Hence $T''$ includes all small types of the
form $\pi(u_1,v_1)\pi(u_2,v_2)\tau.$ Using the first assertion of Lemma, we
obtain that $T''$ consists of types
     $\pi(u_1 u_2,v_1 v_2)\tau$ where $\tau\in T.$ Putting $v=v_1 v_2$
     we obtain: $T''=\{\pi(u_1 u_2,v)\tau\mid \tau\in T,l(u_1 u_2)=l(v)\}.$
This set is (again by definition) $\Pi(u_1 u_2)T.$
So Lemma is completely proved.

We proceed to return to the assertion formulated at the beginning of this
Section and to prove it.

{\bf Theorem 1.}
{\it
 $L\subset\Sigma_1^*$, the set of words-solutions for the equation
(\ref{ur_otn_u}), is a regular language.}

{\bf Proof.} Consider the following finite automaton $G$. Its vertices are
various large types. An arrow goes from T to $T'$ and is marked by the
symbol $x$ iff $T'=\Pi(x)T.$ The initial vertex is
   $T(\lambda),$ and final vertices are various good large types. Some
$u=x_k\cdots x_0$ is a solution iff
   $T(u)$ is a good type. Hence $\Pi(u)T(\lambda)=
   \Pi(x_k)\cdots\Pi(x_0)T(\lambda)$ is a good type, and this in turn is
equivalent to the assertion that $T(u)$ is a final vertex and the end of the
path $x_k\cdots x_0.$
   Thus $u$ is a solution iff $u$ belongs to the language represented by the finite automaton $G.$ This implies the assertion of Lemma.

\section{ Equations over a matrix ring}

Any element of the ring ${\bf M_n({\bf R})}$ nay be interpreted in two ways: as a polynomial in  $\vartheta_1,\ldots,\vartheta_r$ with coefficients from
${\bf M_n({\bf Z_p})}$, and as a matrix with entries from ${\bf R}.$ Suppose $f(B)$ is a polynomial in a matrix $B$ with coefficients from
${\bf R}$ (the matrix itself belongs to ${\bf M_n({\bf R})}$).
Denote the ring consisting of such polynomials by ${\bf R}[B].$
Let $ deg$~$B,$ $B\in
{\bf M_n({\bf R})},$ be the sum of powers of $B$ as a polynomial in
$\vartheta_1,\ldots,\vartheta_r$ with matrix coefficients.

{\bf Definition.}
 A matrix $B\in {\bf M_n({\bf F})}$ is {\it rational of standard form} if it
has the form\\
$$
  \left(
  \begin{array}{ccccc}
  0 & 0 &\cdots & 0 & f_0\\
  1 & 0 &\cdots & 0 & f_1\\
  0 & 1 &\cdots & 0 & f_2\\
  \vdots &\vdots &\ddots &\vdots &\vdots\\
  0 & 0 &\cdots & 1 &f_{n-1}
  \end{array}
  \right)
$$\\
where the polynomial $\xi^n-\sum\limits_{i=0}^{n-1} f_i \xi^i$
is irreducible and separable over {\bf F}.\\
In fact this the matrix of multiplication by $\xi$ in the algebraic
extension of the field {\bf F} by a root of the above polynomial in the
basis of extension consisting of powers of this root
(see [2, бва.~429-455]).

{\bf Definition.}
 A matrix $B'\in {\bf M_n({\bf R})}$ is {\it entire of standard form} if it
has the form\\
$$
  \left(
  \begin{array}{ccccc}
  0 & 0 &\cdots & 0 &\varrho_0\\
  \varrho & 0 &\cdots & 0 &\varrho_1\\
  0 &\varrho &\cdots & 0 &\varrho_2\\
  \vdots &\vdots &\ddots &\vdots &\vdots\\
  0 & 0 &\cdots &\varrho &\varrho_{n-1}
  \end{array}
  \right)
$$\\
 where the polynomial ${\xi}^n-\sum\limits_{i=0}^{n-1}
\frac{\varrho_i}{\varrho} \xi^i$
is irreducible and separable over {\bf F}.\\

{\bf Remark.}
 If a matrix $B'$ is entire of standard form then there exists a unique
matrix $B$, rational of standard form, such that $B'=\varrho B,\varrho\in
{\bf R}.$

Let ${\bf F}[B]$ be the ring of polynomials having the following form:\\
$\sum\limits_{k=0}^{m} f_k(\overline\vartheta)B^{k}(\overline\vartheta),
 f_k\in {\bf F} , B\in{\bf M_n({\bf F})}.$ We have the following

{\bf Lemma 7 (on simplifying the form of an equation).}\\
{\it The following assertions are equivalent:\\
 a) The set of solutions for an EDE over ${\bf M_n({\bf A})}$ is a regular
language.\\
 b) The set of solutions for an EDE over {\bf A} is a regular language.\\
 c) The set of solutions for an EDE over a finite algebraic extension
{\bf F} is a regular language.\\
 d) The set of solutions for an EDE over a finite separable algebraic
extension
{\bf F} is a regular language.\\
 e) The set of solutions for an EDE over ${\bf F}[B]$ where $B$ is a
rational matrix of standard form is a regular language.\\
 f) The set of solutions for an EDE over ${\bf R}[B']$ where $B'$ is an
entire matrix of standard form is a regular language.
}

{\bf Proof.}
  Observe obvious implications:
  a) $\Rightarrow$ b) $\Rightarrow$ c) $\Rightarrow$ d).\\
We proceed to prove c) $\Rightarrow$ b).
  Consider an EDE over {\bf A}. It involves only a finite number of
coefficients, and all of them are algebraic over {\bf F}. Hence they belong
to a finite algebraic extension of {\bf F}, and the original equation is an
EDE over this extension.\\

  Now we shall prove d) $\Rightarrow$ c). Consider an EDE over a finite
algebraic extension of {\bf F}:
  $$
    \sum_{i=1}^{s} b_i a_{i1}^{n_1} \cdots a_{it}^{n_t}=0.
  $$
For some $M$ all of $b_i^{p^M}, a_{ik}^{p^M}$ are separable over {\bf F}.
Consider the equation
  $$
    \sum_{i=1}^{s} (b_i)^{p^M} (a_{i1})^{p^M n_1} \cdots (a_{it})^{p^M n_t}=
    (\sum_{i=1}^{s} b_i a_{i1}^{n_1} \cdots a_{it}^{n_t})^{p^M}=0.
  $$
This equation is equivalent to the original one and is an EDE over a finite
separable algebraic extension of {\bf F}.\\

     Now we postpone the proof for the most difficult implication
b) $\Rightarrow$ a) and proceed to show that d) and e) are equivalent.\\
  d) $\Rightarrow$ e).
  Since $B$ is rational of standerd form, ${\bf F}[B]$ is a finite separable
algebraic extension for {\bf F}.\\
  e) $\Rightarrow$ d).
By the primitive element theorem, every finite separable algebraic extension
may be represented in the form ${\bf F}[\xi]$
  where $\xi$ is a root of some separable over {\bf F} polynomial
  $z^n-\sum\limits_{i=0}^{n-1} f_i z^i.$ Consider the matrix $B$:\\
    $$
      \left(
      \begin{array}{ccccc}
      0 & 0 &\cdots & 0 & f_0\\
      1 & 0 &\cdots & 0 & f_1\\
      0 & 1 &\cdots & 0 & f_2\\
      \vdots &\vdots &\ddots &\vdots &\vdots\\
      0 & 0 &\cdots & 1 &f_{n-1}
      \end{array}
      \right).
    $$\\
It is rational of standard form, and ${\bf F}[B]$
  is isomorphic to ${\bf F}[\xi].$ Hence every equation from d) is an
equation from e).\\
  We proceed to prove that e) and f) are equivalent.\\
  e) $\Rightarrow$ f).
Consider an EDE over ${\bf R}[B']$:\\
  $$
    \sum_{i=1}^{s} Q_i(B')[{\overline{P_i}}^{\phi(u)}](B')=0.
  $$
  Suppose $f$ is an arbitrary polynomial over {\bf R}.
  Since $B'=\varrho B, \varrho\in R,$ $B$ is rational of standard form, we
have $f(B')=f(\varrho B)\in F[B]),$ hence ${\bf R}[B']\subset {\bf F}[B].$
Thus our equation is an EDE over ${\bf F}[B].$\\
  f) $\Rightarrow$ e)
Consider an EDE over ${\bf F}[B]:$\\
  $$
    \sum_{i=1}^{s} Q_i(B)[{\overline{P_i}}^{\phi(u)}](B)=0.
  $$
Find the common denominator for all $Q_i,P_{ik}$ and put
  $Q_i=\frac{Q'_i}{\sigma}, P_{ik}=\frac{P'_{ik}}{\sigma}$ where
$Q_i,P_{ik},
  \sigma$ are polynomials over {\bf R}.
  We have
  $$
    \left(\frac{1}{\sigma^{1+\phi^{(1)}(u)+\cdots+\phi^{(t)}(u)}}\right)
    \sum_{i=1}^{s} Q'_i(B)[{\overline{P'_i}}^{\phi(u)}](B)=0.
  $$
This equation is equivalent to the following:
  $$
    \sum_{i=1}^{s} Q'_i(B)[{\overline{P'_i}}^{\phi(u)}](B)=0.
  $$
Now use the fact that $B=\frac{1}{\varrho} B'.$
Hence for any polynomial $f$ over {\bf R} we have
  $f(B)=\frac{f'(B')}{\varrho^m}$ where $f'$ also is a polynomial over {\bf
R},
  $m= deg \mbox{ } f.$
Then the above equation may be written in the form\\
  $$
    \left(\frac{1}{\varrho^{m(1+\phi^{(1)}(u)+\cdots+\phi^{(i)}(u))}}\right)
    \sum_{i=1}^{s} Q''_i(B')[{\overline{P''_i}}^{\phi(u)}](B')=0.
  $$
Multipying by $\varrho^{m(1+\phi^{(1)}(u)+\cdots+\phi^{(i)}(u))},$
  we obtain an EDE over ${\bf R}[B']$ equivalent to the original one.\\
  Finally we prove the last implication.\\
  b) $\Rightarrow$ a). Consider an EDE over $M_n({\bf A})$:\\
    $$
      \sum_{i=1}^{s} B_{i0}A_{i1}^{n_1}B_{i1}\cdots A_{it}^{n_t}B_{it}=0
    $$
   where $ A_{ik},B_{il}\in {\bf M_n({\bf A})}.$
   Since {\bf A} is algebraically closed,
   $A_{ik}$ is representable in the form $C_{ik}A_{\mbox{J}ik}C_{ik}^{-1}$
   where $A_{\mbox{J}ik}$ is a Jordan matrix.
   Then $A_{\mbox{J}ik}=D_{ik}+R_{ik}$ where $D_{ik}$ is a diagonal matrix
and
   $R_{ik}$ is nilpotent. Hence there exists
   $M$ such that for any $i, k : R_{ik}^{p^M}=0$
   and so $A_{\mbox{J}ik}^{p^M}=(D_{ik}+R_{ik})^{p^M}=D_{ik}^{p^M}.$\\
   Represent all $n_i$ in the form $n'_i+p^M n_i^*$ where all $n'_i\le p^M.$
   Denote $D_{ik}^{p^M}$ by ${D'}_{ik}.$ In the new notation, the original
equation takeas the form\\
   $$
     \sum_{i=1}^{s} (B_{i0}C_{i1}A_{\mbox{J}i1}^{n'_1}) {D'}_{i1}^{n_1^{*}}
(C_{i1}^{-1}B_{i1}C_{i2}A_{\mbox{J}i2}^{n'_2})
{D'}_{i2}^{n_2^{*}}\cdots {D'}_{it}^{n_t^{*}} (C_{it}^{-1}B_{it})=0
   $$
   or for fixed $n'_1,\ldots,n'_t$:\\
   $$
     \sum_{i=1}^{s} {B'}_{i0}{D'}_{i1}^{n_1^{*}}{B'}_{i1}
     \cdots {D'}_{it}^{n_t^{*}}{B'}_{it}=0.
   $$
   Denote the $kl$-th entry of ${B'}_{ij}$ by $\beta_{ij,kl}\in {\bf A},$
   and $k$-th entry of the diagonal matrix ${D'}_{ij}$ by
   $\lambda_{ij,k}\in {\bf A}.$
   Let $\sigma(x,y,\overline{n^*})$ denote the expression
   $ \sum\limits_{1\le z_1,\ldots,z_t\le n,1\le i\le s}
     \beta_{i0,xz_1}\beta_{i1,z_1 z_2}\cdots\beta_{it,z_t y}
     \lambda_{i1,z_1}^{n_1^{*}}\cdots\lambda_{it,z_t}^{n_t^{*}}$
   for various values of $x,y$ from $1$ to $n.$
   Then the above equation is equivalent to the following system
   (depending on $\overline n'- \langle n'_{1},\ldots,n'_{n}\rangle $)~:
   $$
     \sigma(x,y,\overline{n^*})=0.
   $$
   This is a system of $n^2$ EDE over {\bf A} and so by b)
   $$
   L(\overline{n'},x,y)=\{u\in\Sigma_1^*\mid\sigma (x,y,\phi(u))=0 \}
   $$
   is a regular language. Then
   $$
   L(\overline{n'})=
   \bigcap\limits_{1\le x,y\le n} L(\overline{n'},x,y)
   $$
also is a regular language, and hence the set of solutions for the original
EDE
   $$
   L=\bigcup\limits_{0\le n'_1,\ldots,n'_t <  p^M}
   \{\overline{n'}\}*L(\overline{n'})
   $$
is a regular language (here we apply Cleenee's theorem).\\
Lemma is completely proved.\\

Arguing for the ring of polynomials over a field, we have essentially
applied the identity
$\{f(\overline\vartheta)\}^{p}=f(\overline\vartheta^{p}).$
For the ring of polynomials over a non-commutative ring, in particular over
a matrix ring, this identity fails. But it turns that our constructions can
by extended to the case of the ring ${\bf R}[B]$
considered in assertion f) of the above lemma, by means of the following
statement:\\

{\bf Lemma 8 (on conjugation).}\\
{\it
a) Suppose $B(\overline\vartheta)$ is a rational matrix of standard form.
   Then $B^{p}(\overline\vartheta)=$\\
 $C(\overline\vartheta)B({\overline\vartheta}^p)C^{-1}(\overline\vartheta)$
   where $C\in {\bf M_n({\bf F})}.$\\
b) Suppose $B'(\overline\vartheta)$ is an entire matrix of standard form.
   Then ${B'}^{p}(\overline\vartheta)=$\\
   $C'(\overline\vartheta)B'({\overline\vartheta}^p)
{C'}^{-1}(\overline\vartheta)$
   where $C'\in {\bf M_n({\bf R})}.$
}

{\bf Proof.}

\nobreak
  )
      $$
             B=
             \left(
             \begin{array}{ccccc}
             0 & 0 &\cdots & 0 & f_0\\
             1 & 0 &\cdots & 0 & f_1\\
             0 & 1 &\cdots & 0 & f_2\\
             \vdots &\vdots &\ddots &\vdots &\vdots\\
             0 & 0 &\cdots & 1 &f_{n-1}
             \end{array}
             \right).
      $$\\
     Suppose
     $\xi^n-\sum\limits_{i=0}^{n-1} f_i\xi^i=0.$
     Then $B(\overline\vartheta)$ is the matrix of the operator
$x\mapsto\xi x$
     in the basis $1,\xi,\ldots,\xi^{n-1} .$ The matrix
$B^p(\overline\vartheta)$ corresponds to the operator $x\mapsto\xi^p x$ in
the same basis. The matrix $B(\overline\vartheta^p)$ corresponds to the
operator $x\mapsto\xi^p x$ in the basis $1,\xi^p,\ldots,\xi^{p(n-1)}.$
     Due to well-known theorem of linear algebra, these matrices are
conjugate.

  b) The matrix $B'$ is entire of standard form. It is known that
      $B'=\varrho B$ where $B$ is rational of standard form,
      $\varrho\in R.$ Then
      $B'({\overline\vartheta}^p)=\varrho^p B^p(\overline\vartheta)=
      \varrho^p C(\overline\vartheta)B({\overline\vartheta}^p)
      C^{-1}(\overline\vartheta)$
      $C(\overline\vartheta)=\frac{1}{\sigma(\overline\vartheta)}
      C'(\overline\vartheta), C'\in {\bf M_n({\bf R})},
      \sigma(\overline\vartheta)$ is the greatest common divisor for
denominators of all entries of $C(\overline\vartheta).$
      Then $C^{-1}(\overline\vartheta)$=
      $\sigma(\overline\vartheta){C'}^{-1}(\overline\vartheta),$
      and so
      $B'{}^p(\overline\vartheta)$=$\varrho^p
      C'(\overline\vartheta)B({\overline\vartheta}^p)
      {C'}^{-1}(\overline\vartheta)$=$C'(\overline\vartheta)
      B({\overline\vartheta}){C'}^{-1}(\overline\vartheta).$

Consider an EDE over ${\bf R}[B]$ where $B$ is an entire matrix of standard
form:

\begin{equation}
  \label{edu2}
  \sum_{i=1}^{s} Q_i(B)[{\overline{P_i}}^{\phi(u)}](B)=0.
\end{equation}

{\bf Remark.}
 Suppose $f(\xi)$ is an arbitrary polynomial from ${\bf R}[\xi].$ Then
$$
  f^{p}(B(\overline\vartheta))=
C(\overline\vartheta)f(B({\overline\vartheta}^p))C^{-1}(\overline\vartheta).
$$

Now put $u=x_k\cdots x_0.$\\
Transforming the right side of (\ref{edu2}), we subsequently have
$$
      0=\sum\limits_{i=1}^{s}
Q_i(B(\overline\vartheta))[{\overline{P_i}}^{\phi(x_0)}]
(B(\overline\vartheta))[{\overline{P_i}}^{\phi(x_1)p})](B(\overline\vartheta
))\cdots [{\overline{P_i}}^{\phi(x_k)p^k}](B(\overline\vartheta))=
$$
$$
       =\sum\limits_{i=1}^{s}
Q_i(B(\overline\vartheta))[{\overline{P_i}}^{\phi(x_0)}]
(B(\overline\vartheta))C(\overline\vartheta)[{\overline{P_i}}^{\phi(x_1)}]
        (B({\overline\vartheta}^p))C^{-1}(\overline\vartheta)
        \cdots ( C(\overline\vartheta)C({\overline\vartheta}^p)\cdots
$$
$$
        \cdots C({\overline\vartheta}^{p^k}) )
[{\overline{P_i}}^{\phi(x_k)}]
        (B({\overline\vartheta}^{p^k})) (
C^{-1}({\overline\vartheta}^{p^k})\cdots
        C^{-1}(\overline\vartheta) )=
$$
$$
       =\sum\limits_{i=1}^{s}
Q_i(B(\overline\vartheta))[{\overline{P_i}}^{\phi(x_0)}]
        (B(\overline\vartheta))C(\overline\vartheta)
        \cdots [{\overline{P_i}}^{\phi(x_k)}](B({\overline\vartheta}^{p^k}))
        C({\overline\vartheta}^{p^k})
        C^{-1}({\overline\vartheta}^{p^k})\cdots C^{-1}(\overline\vartheta)
$$ \\
Now multiply the expression on the right side by an invertible matrix\\
$ C(\overline\vartheta)C({\overline\vartheta}^p)\cdots
C({\overline\vartheta}^{p^k}) .$
The resulting equation is equivalent to the original one:
\begin{equation}
  \label{equ}
 \sum\limits_{i=1}^{s} Q_i(B)(
[{\overline{P_i}}^{\phi(x_0)}](B(\overline\vartheta))C(\overline\vartheta)
)\cdots
 ([{\overline{P_i}}^{\phi(x_k)}](B({\overline\vartheta}^{p^k}))C
({\overline\vartheta}^{p^k}) )=0.
\end{equation}

We proceed to generalize constructions of the first step to the matrix
case.\\
Suppose $f\in{\bf M_n({\bf R})}.$
Then as before $f(\overline\vartheta) = \sum\limits_{y\in\Sigma_2}
f_y({\overline\vartheta}^p){\overline\vartheta}^{\psi(y)}.$\\

{\bf Definitions.}\\
a){\it The weeding} by a symbol $y\in\Sigma_2$ is
$\varepsilon_y(f)=f_y(\overline\vartheta).$\\
b){\it The weeding} by a word $v=y_k\cdots y_0$ is
$\varepsilon_v(f)=\varepsilon_{y_k}(\cdots(\varepsilon_{y_0}(f))\cdots)$,
that is, the composition of weedings by letters of the word.\\

{\bf Lemma 9 (properties of the weeding operator).}\\
a) $\varepsilon_y(f+g) = \varepsilon_y(f)+\varepsilon_y(g).$\\
         b) $ deg \mbox{ }\varepsilon_y(f) \le\frac{1}{p}  deg \mbox{ }
f.$\\
         c) $\varepsilon_y( f(\overline\vartheta)g({\overline\vartheta}^p) )
=
            \varepsilon_y(f(\overline\vartheta))g(\overline\vartheta).$ \\
         d) If $c$ is a constant then $f=0$ iff for any
            $v\in\Sigma_2^*$ having length $c$ we have $\varepsilon_v(f)=0
.$\\
The proofs of these properties are similar to those given at the first
step.\\
Special operators are defined for the matrix case as follows:

{\bf Definition.}\\
a) Suppose $x\in\Sigma_1, y\in\Sigma_2.$ {\it The special operator}
$S_{x,y}^{(i)}(f)=\varepsilon_y( f[{\overline{P_i}}^{\phi(x)}](B)C )$
where $C$ is the matrix from Lemma 8.\\
b) Suppose $ u=x_k\cdots x_0, v=y_k\cdots y_0$ are words of equal length
from  $\Sigma_1^*$ and $\Sigma_2^*$ accordingly.
 Then $S_{u,v}^{(i)}(f)=S_{x_k,y_k}^{(i)}\cdots S_{x_0,y_0}^{(i)}(f).$\\

{\bf Lemma 10 (on decreasing the degree).}\\
{\it There exists
 $ N_0$ such that for any $ N'\ge N_0,$ $1\le i\le s,$
 $x\in\Sigma_1,$ $y\in\Sigma_2 $ we have
 $  deg \mbox{ } f\le N'\Rightarrow
 deg \mbox{ } S_{x,y}^{(i)}(f)\le N'.$
In other words, rather high degrees of polynomials can be only decreased by
special operators.}\\
{\bf Proof.}\\
  Denote $max \mbox{ }deg \mbox{ } P_{ik}$ by $M,$
       and $\frac{prM+ deg \mbox{ } C}{p-1}$
       by $N_1.$ Then $N_1$ is the desired $N_0.$
        We proceed to prove this. Suppose
       $ N'=N_1+K.$ Then $deg \mbox{ } f\le N'$ implies
       $$
         deg \mbox{ }f{\overline{P_i}}^{\phi(x)}C\le N'+Mpr+ deg \mbox{ }C=
         (prM+ deg \mbox{ } C)\left(\frac{p}{p-1}\right)+K.
       $$
       Furthermore
       $$
         deg \mbox{ } S_{x,y}^{(i)}(f)\le
       \frac{1}{p}  deg \mbox{ } f{\overline{P_i}}^{\phi(x)}C\le
       \frac{prM+ deg \mbox{ } C}{p-1}+\frac{K}{p}\le N'.
       $$
       Lemma is proved.

{\bf Definitions.}\\
a){\it A small type} $T=(f_1,\ldots,f_s)$ is a string of matrices from ${\bf
M_n({\bf R})}$ such that $ deg \mbox{ } f_i\le N_2 ,
N_2= max \{  max$ $deg \mbox{ }(Q_i),N_0\}.$\\
b)Let $x,y$ be symbols from alphabets $\Sigma_1$ and $\Sigma_2$ accordingly.
{\it The extension} $\pi(x,y)\tau$ of a small type
$\tau=(f_1,\ldots,f_s)$ by these symbols is the small type $\tau'= (f'_1
,\ldots, f'_s),
f'_i=S_{x,y}^{(i)}(f_i).$\\
c)Suppose $u=x_k\cdots x_0, v=y_k\cdots y_0$ are words of length from
$\Sigma_1^*$ and $\Sigma_2^*$ accordingly.
Then {\it the extension} $\pi(u,v)\tau$ of a small type $\tau$ by this pair
of words is the composition of its extensions by pairs of symbols
$\pi(x_k,y_k)\cdots\pi(x_0,y_0)\tau.$\\

{\bf Remark.}
If $\tau=(f_1,\ldots,f_s)$ then $\pi(u,v)\tau=(f'_1,\ldots,f'_s)$
 where $f'_i=S_{u,v}^{i}(f_i).$
This follows immediately from definitions of $\pi(u,v)$ and $S_{u,v}^{i}.$

{\bf Definition.}\\
a){\it A large type} T is an arbitrary set of small types.\\
b) Suppose $u\in\Sigma_1^*.$ {\it The extension} of a type T by the word $u$
is the large type ${\Pi}(u)T=\{\pi(u,v)
\tau \mid \tau\in T,
 v\in\Sigma_2^*, l(v)=l(u)\}.$\\

{\bf Lemma 11.}
{\it
${\Pi}(u_1 u_2)T={\Pi}(u_1){\Pi}(u_2)T,$ that is, an extension of a large
type by a concatenation of two words is composition of extensions of this
type by the given words.}\\
{\bf Proof.}
 Denote ${\Pi}(u_2)T$ by $T',$ and ${\Pi}(u_1)T'$ by $T''.$
By definition of extension of a large type,
${T'}=\{\pi(u_2,v_2)\tau \mid \tau\in T, l(v_2)=l(u_2)\}.$\\
In turn, $T''$ is (again by definition)
$\{\pi(u_1,v_1)\tau'\mid \tau'\in T', l(v_1)=l(u_1)\}$
and, by virtue of formula for $T',$ equals
$\{\pi(u_1,v_1)\pi(u_2,v_2)\tau \mid \tau\in T,
 l(v_1)=l(u_1),$ $l(v_2)=l(u_2)\},$\\ or $\{\pi(u_1 u_2,v)\tau \mid \tau\in
T,
l(v)=l(u_1 u_2)\}.$\\ This precisely coincides with ${\Pi}(u_1 u_2)T.$
 Lemma is proved.\\

{\bf Definition.}
 Let $u,v$ be words of equal length from $\Sigma_1^*$ and $\Sigma_2^*$
accordingly.\\
 a){\it The small type of the pair of words}
 $\tau(u,v)$ is $\pi(u,v)\tau(\lambda,\lambda)$ where
 $\tau(\lambda,\lambda)=(Q_1(B),\ldots,Q_s(B)).$\\
 b){\it The large type of the word} $T(u)$ is $\{\tau(u,w) \mid
l(w)=l(u)\}.$

{\bf Lemma 12.}
{\it
 $T(u)={\Pi}(u)T(\lambda),$ that is, the large type of the word $u$
is the extension by this word of the large type of the empty word.}\\
{\bf Proof.}\\
By definition, $T(\lambda)=\{\tau(\lambda,\lambda)\}.$
So, using only definitions for extensions of large and small types, we
easily obtain
 $$
   {\Pi}(u)T(\lambda) = \{\pi(u,v)\tau(\lambda,\lambda)\mid l(v)=l(u)\}
   = \{\tau(u,v)\mid l(v)=l(u)\} = T(u).
 $$
Lemma is proved.\\

{\bf Definition.}
a) A small type $\tau=(f_1,\ldots,f_s)$ is {\it good} if
 $\sum\limits_{i=1}^{s} f_i=0.$\\
b) A large type T is {\it good} if all $\tau\in T$ are good.\\

{\bf Theorem 3.}
{\it
 Suppose $u$ is an arbitrary word from $\Sigma_1^*.$  Then $T(u)$ is a good
type iff $u$ is a solution for the original EDE.}

{\bf Proof.}
Define the length of $u$ by $c.$ By definition of $T(u)$, it is a good type
if for any $v\in\Sigma_2^*$ of length $c$ the type
$\tau(u,v)=\pi(u,v)\tau(\lambda,\lambda)$ is good. This means in turn that
for any $v\in\Sigma_2^*$ of length $c$ we have
  $$
    \sum\limits_{i=1}^{s} S_{u,v}^{(i)}( Q_i(B(\overline\vartheta)) )=0.
  $$
This implies that for all $v\in\Sigma_2^*$ of length $c$ we have
  $$
    \sum\limits_{i=1}^{s} \varepsilon_v(
Q_i(B)[{\overline{P_i}}^{\phi(x_0)}]
    (B(\overline\vartheta))C(\overline\vartheta)
    \cdots [{\overline{P_i}}^{\phi(x_c)}](B({\overline\vartheta}^{p^c}))
    C({\overline\vartheta}^{p^c}) )=0,
  $$
or for any $v\in\Sigma_2^*$ of length $c$
 $$
    \varepsilon_v( \sum\limits_{i=1}^{s}
Q_i(B)[{\overline{P_i}}^{\phi(x_0)}]
    (B(\overline\vartheta))C(\overline\vartheta)\cdots
    [{\overline{P_i}}^{\phi(x_c)}](B({\overline\vartheta}^{p^c}))
    C({\overline\vartheta}^{p^c}) )=0.
  $$
Using ae property of weeding (assertion {\it d} of Lemma 9 ), we have now
  $$
  \sum\limits_{i=1}^{s} Q_i(B)[{\overline{P_i}}^{\phi(x_0)}]
  (B(\overline\vartheta))C(\overline\vartheta)\cdots
  [{\overline{P_i}}^{\phi(x_c)}](B({\overline\vartheta}^{p^c}))
  C({\overline\vartheta}^{p^c}) =0,
  $$
  that is, $u$ is a solution of (\ref{equ}), and so $u$ is a solution of
(\ref{edu2}). Theorem is proved.\\
Now it remains to construct the desired fiite automaton.

{\bf Remark.}
The number of large and small types is finite. Small types are matrices of
order $n,$ their entries are polynomials of bounded degree in
$r$ variables over ${\bf Z_p.}$ Large types are subsets of some finite
sets.\\

{\bf Theorem 4.}
{\it The set of solutions for an EDE over the ring ${\bf R}[B]$ where B is
an entire matrix of standard form, is a regular language.}\\
Construction of the finite automaton is completely similar to the case of
the ring of polynomials. Again vertices are large types, and an arrow marked
by $x$ goes from $T_1$ to $T_2$ iff $T_2=\Pi(x)T_1.$ The initial vertex is
$T(\lambda),$ and final vertices are all of good large types.
The proof is similar to the one given in the preceding section.

Theorem 4 and Lemma 7 immediately imply the following

{\bf Theorem 5.}
{\it
 Suppose {\bf F} is a field, char~${\bf F}=p,$
then the set of solutions for an EDE over ${\bf M_n({\bf F})}$ is a regular
language.}

{\bf Proof.}
An EDE includes a finite number of matrix entries, so all of them belong to
some finite extension of ${\bf Z_p.}$ Any such extension may be included
in {\bf A} if $r$ is its transcendence degree. Hence the original equation is an EDE over ${\bf M_n({\bf A})}.$  By Lemma 7 and Theorem 4 we obtain that the set of its solutions is a regular language.

{\bf Corollary.}
{\it
 If {\bf R} is a ring representable by matrices over a field
{\bf F}, $ char$
$ {\bf F}=p,$ then the set of solutions for an EDE over {\bf R} is a regular
language.}
\\
\\

\end{document}